\documentclass{amsart}
\usepackage{amssymb}
\usepackage{braket}
\usepackage{mathrsfs}
\usepackage{ifthen}

\usepackage{tikz}
\usetikzlibrary{patterns}
\usepackage{comment} 
\usepackage{mleftright}
\usepackage{hyperref}
\usepackage{cleveref}
 \usepackage[all]{xy}
 \usepackage{amscd}
 \usepackage{amsrefs}
 \usepackage{color}
 \usepackage{enumitem}
\newlist{steps}{enumerate}{1}
\setlist[steps, 1]{label = Step \arabic*:}
\usepackage[abs]{overpic}
\usepackage[latin1]{inputenc}

\makeatletter
\DeclareRobustCommand\widecheck[1]{{\mathpalette\@widecheck{#1}}}
\def\@widecheck#1#2{%
   \setbox\z@\hbox{\m@th$#1#2$}%
   \setbox\tw@\hbox{\m@th$#1%
      \widehat{%
         \vrule\@width\z@\@height\ht\z@
         \vrule\@height\z@\@width\wd\z@}$}%
   \dp\tw@-\ht\z@
   \@tempdima\ht\z@ \advance\@tempdima2\ht\tw@ \divide\@tempdima\thr@@
   \setbox\tw@\hbox{%
      \raise\@tempdima\hbox{\scalebox{1}[-1]{\lower\@tempdima\box\tw@}}}%
   {\ooalign{\box\tw@ \cr \box\z@}}}
\makeatother

\theoremstyle{plain}
\newtheorem{thm}{Theorem}[section]
\crefname{thm}{Theorem}{Theorems}
\Crefname{thm}{Theorem}{Theorems}

\crefname{prop}{Proposition}{Propositions}
\Crefname{prop}{Proposition}{Propositions}

\crefname{lem}{Lemma}{Lemmas}
\Crefname{lem}{Lemma}{Lemmas}
\newtheorem{cor}[thm]{Corollary}
\crefname{cor}{Corollary}{Corollaries}
\Crefname{cor}{Corollary}{Corollaries}

\crefname{claim}{Claim}{Claims}
\Crefname{claim}{Claim}{Claims}

\crefname{property}{Property}{Properties}
\Crefname{property}{Property}{Properties}

\crefname{problem}{Problem}{Problems}
\Crefname{problem}{Problem}{Problems}

\crefname{conjecture}{Conjecture}{Conjecture}
\Crefname{conjecture}{Conjecture}{Conjecture}

\theoremstyle{definition}

\crefname{defn}{Definition}{Definitions}
\Crefname{defn}{Definition}{Definitions}

\crefname{notation}{Notation}{Notations}
\Crefname{notation}{Notation}{Notations}

\crefname{convention}{Convention}{Conventions}
\Crefname{convention}{Convention}{Conventions}

\crefname{cond}{Condition}{Conditions}
\Crefname{cond}{Condition}{Conditions}

\crefname{assum}{Assumption}{Assumptions}
\Crefname{assum}{Assumption}{Assumptions}

\crefname{conj}{Conjecture}{Conjectures}
\Crefname{conj}{Conjecture}{Conjectures}

\theoremstyle{remark}
\newtheorem{rem}[thm]{Remark}
\crefname{rem}{Remark}{Remarks}
\Crefname{rem}{Remark}{Remarks}

\crefname{ex}{Example}{Examples}
\Crefname{ex}{Example}{Examples}

\crefname{section}{Section}{Sections}
\Crefname{section}{Section}{Sections}
\crefname{subsection}{Subsection}{Subsections}
\Crefname{subsection}{Subsection}{Subsections}
\crefname{figure}{Figure}{Figures}
\Crefname{figure}{Figure}{Figures}

\newtheorem*{acknowledgement}{Acknowledgement}

\newcommand{\Z}{\mathbb{Z}}

\newcommand{\fraks}{\mathfrak{s}}

\newcommand{\inc}{\hookrightarrow}

\newcommand{\id}{\mathrm{id}}

\newcommand{\s}{\mathfrak{s}}

\def\om{\omega}

\def\Spinc{\text{Spin}^c}

\def\det{\operatorname{det}}

\def\id{\operatorname{Id}}

\newcommand{\mbar}[1]{{\ooalign{\hfil#1\hfil\crcr\raise.167ex\hbox{--}}}}

\title[Generalized Thurston-Bennequin inequality]{A note on generalized Thurston--Bennequin inequalities}

\author{Nobuo Iida}
\address{Department of Mathematics Tokyo Institute of Technology 2-12-1, Ookayama, Meguro, Tokyo 152-8551 Japan}
\email{iida.n.ad@m.titech.ac.jp}

\author{Hokuto Konno}
\address{Graduate School of Mathematical Sciences, the University of Tokyo, 3-8-1 Komaba, Meguro, Tokyo 153-8914, Japan}
\email{konno@ms.u-tokyo.ac.jp}

\author{Masaki Taniguchi}
\address{2-1 Hirosawa, Wako, Saitama 351-0198, Japan}
\email{masaki.taniguchi@riken.jp}

\date{\today}
\begin{document}

\maketitle
\begin{abstract}
We give a generalized Thurston--Bennequin-type inequality for links in $S^3$ using a Bauer--Furuta-type invariant for 4-manifolds with contact boundary.
As a special case,
we also give an adjunction inequality for smoothly embedded orientable surfaces with negative intersection in a closed oriented smooth 4-manifold whose non-equivariant Bauer--Furuta invariant is non-zero.
\end{abstract}


\section{Main results}

Adjunction inequalities are lower bounds on genera of smoothly embedded surfaces in a 4-manifold modeled on the adjunction formulae for algebraic curves in an algebraic surface, and they have been important tools to study 4-dimensional topology.
In this paper, we give a new adjunction-type inequality.
The most general result is given as a generalized Thurston--Bennequin-type inequality for links in $S^3$
(\cref{thm: relative main}), and it deduces also an adjunction inequality for closed surfaces with negative self-intersection in a closed 4-manifold, described in \cref{Adjunction inequality for negative-self intersections}.

\subsection{Adjunction inequality for negative-self intersections}
\label{Adjunction inequality for negative-self intersections}

First, we describe our adjuction-type inequality for closed 4-manifolds.

Most known adjunction inequalities are proved for surfaces with non-negative self-intersection (see, for example, \cite{Law97}).
However, Ozsv{\'a}th--Szab{\'o} \cite[Corollary1.7]{OS00} proved the adjuction inequality for surfaces with negative self-intersection provided that the 4-manifold has non-vanishing Seiberg--Witten invariant and is of Seiberg--Witten simple type.
Recently, Kato, Nakamura, and Yasui~\cite[Theorem~1.7]{KNY20} and Baraglia~\cite[Theorem~1.8]{Ba22} proved adjunction inequalities without simple type assumption provided that the non-vanishing of the (mod 2) Seiberg--Witten invariant and $b^+\equiv 3 \operatorname{mod} 4$.
Note that the results by Ozsv{\'a}th--Szab{\'o}, Kato--Nakamura--Yasui, and Baraglia cannot be applied to connected sums of 4-manifolds with $b^+>0$ since the Seiberg--Witten invariants vanish for such 4-manifolds.

We give an adjunction-type inequality for surfaces with negative self-intersection also for some class of 4-manifolds obtained as connected sums of 4-manifolds with $b^+>0$.
Also, our method does not need to require Seiberg--Witten simple type as in \cite{KNY20,Ba22}.
The statement is formulated in terms of the (non-equivariant) Bauer--Furuta invariant \cite{BF04}:

\begin{thm}
\label{thm: main thm closed}
Let $X$ be an oriented closed smooth 4-manifold with $b_1(X)=0$ and $\fraks$ be a Spin$^c$ structure on $X$.
Suppose that the non-equivariant Bauer--Furuta invariant $BF(X,\fraks)$ is non-trivial.
Let $\Sigma \inc X$ be a smoothly embedded oriented closed surface with $g(\Sigma)>0$.
Then  we have
\begin{align}
\label{eq: main adj closed}
|c_1(\fraks) \cdot [\Sigma]| + [\Sigma]^2\leq 2g(\Sigma).
\end{align}
\end{thm}

By Bauer's calculation \cite{B04} of Bauer--Furuta invariants, we may apply \cref{thm: main thm closed} to many concrete 4-manifolds.
Recall that the (Seiberg--Witten) formal dimension $d(\fraks)$ for a closed Spin$^c$ 4-manifold $(X,\fraks)$ is defined as 
\[
d(\fraks) = \frac{c_1(\fraks)^2-2\chi(X)-3\sigma(X)}{4}.
\]

\begin{cor}
\label{main cor closed}
Let $1\leq n \leq 3$.
For $i \in \{1, \ldots,n\}$, let $X_i$ be an oriented closed smooth 4-manifold with $b^+(X_i) \equiv 3 \mathrm{\ mod\ } 4$ and $b_1(X_i)=0$, and $\fraks_i$ be a Spin$^c$ structure on $X_i$ with formal dimension 0.
Suppose that the Seiberg--Witten invariants $SW(X_i,\fraks_i)$ are odd numbers for all $i$.
Set 
\[
(X,\fraks) = \#_{i=1}^n (X_i,\fraks_i)
\]
and let $\Sigma \inc X$ be a smoothly embedded oriented closed surface with $g(\Sigma)>0$.
Then  we have
\[
|c_1(\fraks) \cdot [\Sigma]| + [\Sigma]^2
\leq 2g(\Sigma).
\]
\end{cor}

\begin{rem}
We describe the nature of \cref{thm: main thm closed,main cor closed}.
\begin{enumerate}
    \item In the case that $[\Sigma]^2 \geq 0$, the results of \cref{thm: main thm closed,main cor closed} have been known (see such as \cite{IMT21,Ba22}).
    The new parts of \cref{thm: main thm closed,main cor closed} are the statements for the case that $[\Sigma]^2<0$.
    \item \cref{thm: main thm closed,main cor closed} generalizes the adjuction inequality by Baraglia~\cite[Theorem~1.8]{Ba22} from 4-manifolds with non-trivial (mod 2) Seiberg--Witten invariant to 4-manifolds with  non-trivial Bauer--Furuta invariant.
    \item As in \cite[Theorem~1.8]{Ba22}, the inequality \eqref{eq: main adj closed} is slightly weaker than the usual adjunction inequality: the right-hand side of the usual adjunction inequality is $2g(\Sigma)-2$.
    \item For the case that $n=1$ in \cref{main cor closed}, the statement of the \lcnamecref{main cor closed} follows from a special case of the Ozsv{\'a}th--Szab{\'o}'s adjuction inequality \cite{OS00} for negative-self intersection and the usual adjunction inequality for non-positive self-intersection, {\it provided that} $(X,\fraks)$ is of simple type.
Note that we do not have to impose any simple type assumption in \cref{thm: main thm closed,main cor closed} as in \cite[Theorem~1.7]{KNY20} and  \cite[Theorem~1.8]{Ba22}.
\end{enumerate}
\end{rem}

\subsection{Generalized Thurston--Bennequin-type inequality}
\cref{thm: main thm closed} is generalized to a relative genus bound.
Let $X$ be a compact oriented 4-manifold whose boundary is $S^3$ and let $\xi_{std}$ be the standard contact structure on $\partial X=S^3$.
Let $L=\amalg^n_{i=1} K_i \subset (S^3, \xi_{std})$ be an oriented  Legendrian link and $\Sigma \subset X$ be a smooth connected oriented surface whose oriented boundary is $K$.
We define the {\it Thurston--Bennequin number} of $L$ by $\displaystyle tb(L) := \sum_{i=1}^n (tb(K_i) +  \sum_{j \neq i}\operatorname{lk}(K_i , K_j))$.  If we take a front projection of $L$, we can write $\displaystyle tb(L) = \sum_{i=1}^n \left( w(K_i) - \frac{1}{2} \# \operatorname{cusp}(K_i) + \sum_{j \neq i}\operatorname{lk}(K_i , K_j) \right)$, which is the same as $tb(L)$ in \cite{Ca18}, where $w$ is the writhe and $\operatorname{lk}$ means the linking number. 

Since $H^2(S^3;\Z)=H^1(S^3;\Z)=0$, we have canonical isomorphism between the determinant line bundle $\det{\s}$ of the positive spinor bundle and $\xi_{std}$.
The {\it generalized rotation number } $r(L, [\Sigma], \s)$
is defined to be $\langle c_1(\det{\s}, \dot{L}), [\Sigma] \rangle $, 
where $h$ is the isomorphism   and $\dot{L}$ is a section of  $\det{\s}|_{S^3}=\xi_{std}$ on $L$ given by the tangent vector field of $L$, which is unique up to pointwise positive scaling. 
Note that the relative first Chern class with respect to this section $c_1(\det{\s}, \dot{L})$ is an element of $H^2(X, L; \Z)$.

Now we are ready to state our relative adjunction inequality:

\begin{thm}\label{thm: relative main}
Let $X$ be an oriented closed smooth 4-manifold with $b_1(X)=0$ and $\fraks$ be a Spin$^c$ structure on $X$.
Let $L$ be an oriented Legendrian link in $(S^3, \xi)$, where $\xi$ is the standard contact structure.
Suppose that the non-equivariant Bauer--Furuta invariant $BF(X,\fraks)$ is non-trivial.
Let $\Sigma \inc X\setminus \operatorname{int}D^4$ be a smoothly and properly embedded oriented compact connected surface with $g(\Sigma)>0$ such that $\partial \Sigma =L$. Then, we have 
\begin{align}
\label{eq: main adj relative'}
| r(L, [\Sigma], \s, h) | + [\Sigma]\cdot  [\Sigma]+  tb(L) -n+2 \leq 2g(\Sigma) .  
\end{align}
\end{thm}
\begin{rem}
We describe the nature of \cref{thm: relative main}.
\begin{enumerate}
\item We first note that our inequality typically can be applied to symplectic caps $X$ whose boundary is $(S^3, \xi_{std})$  with $b^+(X)\equiv 3 \operatorname{mod}$ 4 and $b_1(X)=0$.
The usual generalized Thurston--Bennequin inequality (\cite{MR06}) holds for symplectic fillings but not for symplectic caps. 
\item Unlike the closed case (\cref{thm: main thm closed,main cor closed}), the result of \cref{thm: relative main} is new also for non-negative self-intersection surfaces.
As a known relative genus bounds for connected sums, in \cite{IMT21}, Mukherjee and the first and third authors gave an adjunction-type inequality for certain connected sums for surfaces with non-negative self-intersections.
Let us compare \eqref{eq: main adj relative'} with the adjunction inequality given in \cite[Remark 5.4]{IMT21}. 
For example, \cite[Remark 5.4]{IMT21} implies 
\[
 [\Sigma]^2  \leq 2g(\Sigma) -2 
\]
under the same assumption of \cref{thm: relative main} with $TB(K)>0$, $c_1(\fraks)=0$ and $[\Sigma]^2 \geq 0$. On the other hand, \eqref{eq: main adj relative'} implies 
\[
 [\Sigma]^2 +  TB(K) \leq 2g (\Sigma)-1, 
\]
where $TB(K)$ is the {\it maximal Thurston--Bennequin number}. 
Therefore, if $TB(K)>1$, the bound from \eqref{eq: main adj relative'} is stronger than that in \cite{IMT21}.
Also we do not need to assume $[\Sigma]^2 \geq 0$ and $TB(K)>0$. 
\item Define the 4-genus $g_4(L)$ as minimal genus of all connected properly and smoothly embedded surfaces in $D^4$ bounded by $L$.
If $X=D^4$, the inequality \eqref{eq: main adj relative'} recovers
\begin{align}
\label{D4case}
|r(L) | + tb(L)   \leq 2g_4(L) + n -2, 
\end{align}
which is known as the {\it Thurston--Bennequin inquality} (see for example \cite[Section 1.2, page 18]{OS04}). Here $r(L)$ is the rotation number with respect to the standard contact structure on $S^3$. If $K$ is {\it Lagrangian fillable}, it is proven in \cite{Ba10} that the inequality \eqref{D4case} is attained by the equality, and therefore the inequality \eqref{D4case} is the best possible inequality.
Moreover, as shown in \cite[page 1949]{Ca18}, this bound \eqref{D4case} is sharp for a positive torus link.
\end{enumerate}
\end{rem}

For other relative adjunction inequalities, see also \cite{MMP20, HR20}, for example. Corresponding to \cref{main cor closed}, we have:

\begin{cor}
\label{main cor relative}
Let $1\leq n \leq 3$.
For $i \in \{1, \ldots,n\}$, let $X_i$ be an oriented closed smooth 4-manifold with $b^+(X_i) \equiv 3 \mathrm{\ mod\ } 4$ and $b_1(X_i)=0$, and $\fraks_i$ be a Spin$^c$ structure on $X_i$ with formal dimension 0.
Suppose that the Seiberg--Witten invariants $SW(X_i,\fraks_i)$ are odd numbers for all $i$.
Set 
\[
\displaystyle (X,\fraks) = \#_{i=1}^n (X_i,\fraks_i).
\]
Let $L$ be an oriented Legendrian link in $(S^3, \xi)$, where $\xi$ is the standard contact structure, and $\Sigma \inc X\setminus \operatorname{int}D^4$ be a smoothly and properly embedded oriented compact surface with $g(\Sigma)>0$ such that $\partial \Sigma =L$.
Then  we have
\begin{align}
\label{eq: relative cor}
| r(L, [\Sigma], \s) | + [\Sigma]\cdot  [\Sigma]+  tb(L) -n +2 \leq 2g(\Sigma). 
\end{align}
\end{cor}

Next, we note an application of \cref{thm: relative main} to negative-definite 4-manifolds:

\begin{cor}
\label{main cor relative2}
Let $X$ be a closed smooth oriented negateive-definite 4-manifold with $b_1(X)=0$.
Let $L$ be an oriented Legendrian link in $(S^3, \xi)$, where $\xi$ is the standard contact structure, and $\Sigma \inc X\setminus \operatorname{int}D^4$ be a smoothly and properly embedded oriented compact surface with $g(\Sigma)>0$ such that $\partial \Sigma =L$.
Then  we have
\begin{align}
\label{eq: relative cor2}
\max_{\substack{\fraks \in \mathrm{Spin}^c(X),\\ c_1(\fraks)^2=-b_2(X)}}| r(L, [\Sigma], \s) | + [\Sigma]\cdot  [\Sigma]+  tb(L) -n +2 \leq 2g(\Sigma). 
\end{align}
Here the maximum is taken over spin$^c$ structures $\fraks$ of $X$ that satisfy $c_1(\fraks)^2=-b_2(X)$.
\end{cor}

\begin{rem}
We note that \cref{thm: relative main} recovers a special case of the relative adjunction inequality for quasi-positive knots given by Baraglia~\cite[Theorem 1.2]{Ba22}, and generalize the results to connected sums as follows. If $K$ is  Lagrangian fillable, the equality $tb(K) = 2g_4(K) - 1$ is proven in \cite{Ba10}. Thus, for a closed oriented 4-manifold $X$ with non-zero mod $2$ Seiberg--Witten invariant for a spin$^c$ structure $\fraks$ , $b^+(X)\equiv 3 \operatorname{mod} 4$ and $b_1(X)=0$, and for a properly embedded oriented surface $ \Sigma$ in $X \setminus \operatorname{int} D^4$ bounded  by $K$, our inequality \eqref{eq: relative cor} can be written as
\[
|c_1(\fraks) \cdot [\Sigma]| + [\Sigma]^2 + 2 g_4(K) \leq 2g(\Sigma),
\]
which is the same inequality in \cite[Theorem 1.2]{Ba22}. 
Note that a Lagrangian fillable knot is known to be a quasi-positive knot, which is proven in \cite{KJ15}.
Thus \cref{main cor relative} recovers \cite[Theorem 1.2]{Ba22} for Lagrangian fillable knots.
Moreover, we can treat connected sums of such 4-manifolds. 
For example, we can consider $X' := K3 \# K3\# K3$. Then, for a Lagrangian fillable knot $K$ and a properly embedded oriented surface $ \Sigma$ in $X' \setminus \operatorname{int} D^4$ bounded by $K$, we have 
\[
[\Sigma]\cdot  [\Sigma] +2 g_4(K) \leq 2g(\Sigma),
\]
which is a new bound even for $[\Sigma]\cdot  [\Sigma] \geq 0$.

We also note \cite[Theorem 1.9]{Ba22} can be recovered from \cref{main cor relative2} for Lagrangian fillable knots.  
\end{rem}


\section{Proof of the theorems}
In this section, we will prove  \cref{thm: main thm closed} and \cref{thm: relative main}.
Noting \eqref{eq: main adj relative'},
it is straightforward to see that \cref{thm: main thm closed} follows from \cref{thm: relative main} by considering connected sum with the standard Stein 4-ball $D^4$ and setting $K=U$, the unknot with the standard Legendrian representation, considered as the boundary of the surface obtained as the connected sum of an embedded 2-disk whose boundary is $U$ and the given surface.
Thus we give a proof of \cref{thm: relative main} below.
The invariant introduced in \cite{Iida21} is used in the proof, so let us explain it.
This invariant $\Psi(X, \xi, \s)$ is the Bauer-Furuta type homotopical version of Kronheimer-Mrowka's invariant for a 4-manifold with contact boundary.
While Kronheimer-Mrowka's invariant is defined for an arbitarary compact oriented $\Spinc$ 4-manifold $(X, \s)$ equipped with a contact structure $\xi$ and identification $\s|_{\partial X}=\s_\xi$, the invariant $\Psi(X, \xi, \s)$ is now defined under additional constraint $b_3(X)=0$, in particular the boundary $\partial X$ must be connected.
The invariant is an element 
\[
\Psi(X, \xi, \s)\in \pi^{st}_{d(\s, \xi)}(S^0)/\{\pm 1\}
\]
of the non-equivariant $d(\s, \xi):=\langle e(S^+, \Phi_0), [X, \partial]\rangle$-th stable homotopy group of the sphere $S^0$, defined up to sign, where $\Phi_0$ is a section on the boundary of the positive spinor bundle $S^+$  constructed from  the contact structure as in \cite{KM97}:


\begin{proof}[Proof of \cref{thm: relative main}]
We mainly follow an argument in \cite[Subsubsection~4.2.3]{MR06}.
Denote the components of $L$ by $L = \amalg_{i=1}^n K_i$.
Set 
\[
N = 
\begin{cases}
0,\quad &[\Sigma]^2 \geq 0,\\
-[\Sigma]^2,\quad &[\Sigma]^2<0
\end{cases}.
\]
Form the connected sum of the first component  $K_1$ of $L$ with $\#_N T_{2,3}$, where $T_{2,3}$ is a positive torus knot of type $(2,3)$.
Precisely, we consider Legendrian representations of $K_1$ and $\#_N T_{2,3}$, and form the connected sum along cusps as in \cite[Subsubsection~4.2.3]{MR06}.
Let us set $L':= K_1 \# (\#_N T_{2,3}) \amalg \amalg^{n}_{i=2} K_i$.

Note that $T_{2,3}$ bounds a smoothly embedded genus-1 surface in $D^4$.
Form a boundary connected sum of $\Sigma$ with a null-homologous punctured genus $n$-surface bounded by $\#_N T_{2,3}$ in $X\setminus \operatorname{int}D^4$ and denote the resulting surface by $\Sigma'$, which is bounded by $L'$.
Next, we attach $(tb(K'_i)-1)$-framed 2-handles along $K'_i$ to $X \setminus \operatorname{int}D^4$, where $K'_1=K_1 \# (\#_N T_{2,3})$ and $K'_i=K_i$ for $i \geq 2$. We denote the resulting 4-manifold by $W$. Note that $W$ is decomposed as 
\[
W = X \# W_{(tb(K'_1)-1, \dots, tb(K'_n)-1)}(L'),
\]
where $W_{(tb(K'_1)-1, \dots, tb(K'_n)-1)}(L')$
is the Stein filling \cite{yasha91} obtained from $D^4$ by attaching $n$ Weinstein 2-handles  along $L'$. Denote by $\omega$ the Stein structure. Note that $b_3(W)=0$.
Define 
\[
\Sigma'' := \Sigma' \cup (\text{the cores of the 2-handles of }W_{(tb(K'_1)-1, \dots, tb(K'_n)-1)}(L') ) .
\]

Let $\fraks_\omega$ be the canonical Spin$^c$ strcuture, and
form a Spin$^c$-structure $\mathfrak{s}_W$ on $W$ as the connected sum of $\fraks$ with $\fraks_\om$. We also write by $\xi_\om$ the induced contact structure on $\partial (W_{(tb(K'_1)-1, \dots, tb(K'_n)-1)}(L'))$ from $\om$.

Now we claim that the Bauer--Furuta type invariant $\Psi (W, \mathfrak{s}_W, \xi_\om)$ is non-trivial. Indeed, we can use a connected sum formula as follows:
\begin{align*}
    \Psi (W, \mathfrak{s}_W, \xi_\om) = BF(X, \fraks) \wedge  \Psi (W_{(tb(K'_1)-1, \dots, tb(K'_n)-1)}(L'), \mathfrak{s}_\om, \xi_\om),
\end{align*}
where $BF(X, \fraks)$ is the non-equivariant Bauer--Furuta invariant. 
Since the 4-manifold with boundary $(W_{(tb(K'_1)-1, \dots, tb(K'_n)-1)}(L'), \om)$ is a Stein filling, from \cite{Iida21}, we have 
\[
\Psi (W_{(tb(K'_1)-1, \dots, tb(K'_n)-1)}(L'), \mathfrak{s}_\om, \xi_\om) = \id .
\]
Thus we conclude that $\Psi (W, \mathfrak{s}_W, \xi_\om) $ is the same as $BF(X, \fraks)$, which is non-trivial. 

By the construction $\Sigma''$, we have 
\begin{align*}
g(\Sigma'' )&= g (\Sigma ) + N,\\
[\Sigma'']^2&= [\Sigma']^2 + \sum_{i=1}^n (tb(K_i') -1)+2 \sum_{i<j}\operatorname{lk}(K'_i, K'_j)\\
& = [\Sigma]^2 +tb(L') -n \\
& = [\Sigma]^2 +tb(L)+ 2N -n ,\\
\langle c_1(\fraks_W), [\Sigma''] \rangle & = r(K, [\Sigma], \s).
\end{align*}
Furthermore, from the adjunction inequality proven in \cite{IMT21} for surfaces with non-negative self-intersection applied to $\Sigma''$, we obtain
\[
|\langle c_1(\fraks_W), [\Sigma''] \rangle | + [\Sigma'']^2 \leq 2g(\Sigma'' )-2.
\]
This completes the proof of \cref{thm: relative main}.  
\end{proof}

\begin{acknowledgement}
We would like to thank David Baraglia for giving comments for an earlier version of this paper.
In particular, \cref{main cor relative2} was pointed out by Baraglia.
We wish also thank Kouichi Yasui for informing us about his papers related to adjunction inequalities.

Nobuo Iida was supported by JSPS KAKENHI Grant Numbers 19J23048, 22J00407 and the Program for Leading Graduate Schools, MEXT, Japan.
Hokuto Konno was supported by JSPS KAKENHI Grant Numbers 17H06461, 19K23412, and 21K13785.
 Masaki Taniguchi was supported by JSPS KAKENHI Grant Numbers 20K22319, 22K13921 and RIKEN iTHEMS Program.
\end{acknowledgement}

\bibliographystyle{plain}
\bibliography{tex}

\end{document}